\documentclass[12pt]{amsart}

\usepackage{amsmath}
\usepackage{latexsym}
\usepackage{amssymb}
\usepackage{amsfonts}

\include{PDF}
\usepackage{graphics}

\renewcommand{\proof}{\par\noindent{\it Proof.\ \ }}
\def\qed{\ifmmode\square\else\nolinebreak\hfill
$\square$\fi\par\vskip12pt}

\def\ov{\overline} 
\def\l{\langle} \def\r{\rangle} 
 
 \def\ZZ{\mathbb Z} 

\def\BB{{\mathcal B}}

\def\mod{{\sf mod~}} \def\val{{\sf val}}

\def\Aut{{\sf Aut}} \def\Inn{{\sf Inn}}
\def\Out{{\sf Out}}
\def\Cos{{\sf Cos}} 
\def\bfK{{\bf K}}

 \def\soc{{\sf soc}}

\def\Cay{{\sf Cay}} 
\def\D{{\rm D}} 
\def\S{{\rm S}} \def\G{{\rm G}}

\def\C{{\bf C}}\def\N{{\bf N}}\def\Z{{\bf Z}} \def\O{{\bf O}}

\def\Ome{{\it \Omega}}
\def\Ga{{\it \Gamma}} 
 \def\Del{{\it\Delta}}

\def\a{\alpha} \def\b{\beta}  \def\s{\sigma}
\def\t{\tau}  \def\o{\omega}

 \def\GL{{\rm GL}}

\def\PGammaL{{\rm P\Gamma L}} 
\def\Sym{{\rm Sym}}

\def\PSL{{\rm PSL}}
\def\GL{{\rm GL}} 
\def\AGL{{\rm AGL}}
  
\def\Sz{{\rm Sz}}  
 
  \def\D{{\rm D}} \def\G{{\rm G}}

\newtheorem{theorem}{Theorem}[section]%
\newtheorem{lemma}[theorem]{Lemma}%
\newtheorem{question}[theorem]{Question}%

\begin{document}

\title[Cayley Graphs]
{On Finite Subnormal Cayley Graphs}

\thanks{{\it {2000} Mathematics subject classification}: 05C25, 20B05}

\thanks{This work was partially supported by NSFC(No.61771019).}

\author{Shu Jiao Song}
\address{shujiao.song@ytu.edu.cn}
\address{School of Mathematics and Information Sciences \\
Yantai University\\
P. R. China}


\maketitle

\begin{abstract}
In this paper we introduce and study a type of Cayley graph --
subnormal Cayley graph.
We prove that a subnormal 2-arc transitive Cayley graph is a normal
Cayley graph or a normal cover of a complete bipartite graph
$\bfK_{p^d,p^d}$ with $p$ prime.
Then  we obtain a generic method for constructing half-symmetric 
(namely edge transitive but not arc transitive) Cayley graphs.
\end{abstract}

\section{Introduction}

For a finite group $G$ and a subset $S\subset G$, 
the {\it Cayley digraph} $\Ga=\Cay(G,S)$ is the digraph with vertices
being the elements of $G$ such that $x,y\in G$ are adjacent if and
only if $yx^{-1}\in S$. 
If $S=S^{-1}=\{s^{-1}\mid s\in S\}$, then the adjacency is symmetric
and thus $\Cay(G,S)$ may be viewed as an (undirected) graph, that is, 
a {\it Cayley graph}.
Let 
$$\mbox{$\widehat G=\{\widehat g:\ x\mapsto xg$ for all $x\in G\mid g\in
  G\}$.}$$
Then $\widehat G\leqslant\Aut\Ga$, and $\widehat G$ acts regularly on the vertex set
$G$, so $\Ga$ is vertex-transitive.

A Cayley graph $\Ga=\Cay(G,S)$ is called {\it normal} if
$\widehat G$ is normal in $\Aut\Ga$.
The class of normal Cayley graphs have nice properties and
play important role in studying Cayley graphs, see
  \cite{cubic,Feng,s-trans-Cay,ET-Cay,Praeger-Cay,Xu-Cay} and references
  therein. 
However, there are various interesting classes of Cayley graphs which
  are not normal.

Here we generalize the concept of normal Cayley graphs.
A subgroup $X$ of $Y$ is called {\it subnormal} if there exists a
  sequence of subgroups $X_0,X_1,\dots,X_l$ of $Y$ such that 
$X=X_0\lhd X_1\lhd\dots\lhd X_l=Y$; denoted by $X\lhd\lhd Y$.
A Cayley graph $\Ga=\Cay(G,S)$ is called {\it subnormal} if
$\widehat G$ is subnormal in $\Aut\Ga$, and more generally, $\Ga$ is
  called {\it $Y$-subnormal} if $\widehat G$ is subnormal in $Y$, where
  $\widehat G\leqslant Y\leqslant\Aut\Ga$. 
In the case where $\widehat G\lhd\lhd Y$ and $\Ga$ is $Y$-edge transitive,
$Y$-arc transitive or $(Y,2)$-arc transitive, $\Ga$ is called 
{\it subnormal edge transitive}, {\it subnormal arc transitive}, or
{\it subnormal $2$-arc transitive}, respectively.
(A graph $\Ga$ is called {\it $(Y,2)$-arc transitive} if $Y\leqslant\Ga$ is
  transitive on the set of 2-arcs of $\Ga$.) 
This paper initiates to study the class of subnormal Cayley graphs.

Typical examples of subnormal Cayley graphs are generalized orbital 
graphs of quasiprimitive permutation groups of simple diagonal type
or compound diagonal type, refer to \cite{Praeger-qp}.
The class of 2-arc transitive graphs is one of the central objects in 
algebraic graph theory, see \cite{s-trans-Cay,Praeger-qp} for references. 
It is shown in \cite{s-trans-Cay} that there are only finite many
`basic' 2-arc transitive Cayley graphs of given valency which are not
  normal. 
Here we show that almost all subnormal 2-arc transitive Cayley graphs 
are normal.

Let $\Ga$ be a $(Y,2)$-arc transitive graph with vertex set $V$.
Let $N$ be a normal subgroup of $Y$ which has at least three orbits on
$V$.
Let $\BB$ be the set of $N$-orbits on $V$.
The {\it normal quotient} $\Ga_N$ of $\Ga$ induced by $N$ is the
graph with vertex set $\BB$ such that $B,B'\in\BB$ are adjacent if and
only if some vertex $u\in B$ is adjacent in $\Ga$ to some vertex 
$v\in B'$. 
Then $\Ga_N$ is $(Y/N,2)$-arc transitive of valency equal the valency
of $\Ga$, and so $\Ga$ is a {\it normal cover} of $\Ga_N$.

\begin{theorem}\label{subnormal-2-trans}
Let $\Ga=\Cay(G,S)$ be connected and undirected.
Assume that $Y\leqslant\Aut\Ga$ is such that $\widehat G\lhd\lhd Y$ and 
$\Ga$ is $(Y,2)$-arc transitive.
Then either $\widehat G\lhd Y$, or $\widehat G\lhd X\lhd Y$ with $|Y:X|=2$ 
and $\Ga$ is a normal cover of the complete bipartite
graph $\bfK_{p^d,p^d}$, where $p$ is an odd prime.
\end{theorem}

This theorem is proved in Section~3.
The next theorem is a by-product for proving
Theorem~\ref{subnormal-2-trans}, which extends a classical result for 
primitive permutation groups, that is, \cite[Theorem\,3.2C]{DM-book} and 
\cite[Theorems\,18.4 and 18.5]{Wielandt}, to the general transitive
permutation groups.
Some special cases of this result have been obtained and used in the
study of symmetrical graphs, see for example, 
\cite[Lemma\,2.1]{CI-pty} and \cite[Lemma\,2.1]{FLP-Ree}.

Regular edge transitive graphs are divided into three disjoint classes:
{\it symmetric} (arc transitive); 
{\it semi-symmetric} (vertex intransitive);
{\it half-symmetric} (vertex transitive but not arc transitive).
We remark that in the literature, half-symmetric graphs were called 
{\it half-transitive} graphs.
However, since `half-transitive' is a classical concept for
permutation groups and often occurs in the area of transitive graphs,
we would prefer to call them `half-symmetric' instead of 
`half-transitive'.

Constructing and characterizing half-symmetric graphs is 
an active topic in symmetrical graph theory which has received
considerable attention, see for example \cite{ET-Cay, Li-Song,Marusic, Song}.
The following result provides a generic method for constructing
half-symmetric graphs as subnormal Cayley graphs, which is proved
in Section~4.

\begin{theorem}\label{half-symmetric}
Let $T$ be a finite simple group containing an element $t$ which is
not conjugate in $\Aut(T)$ to $t^{-1}$.
Let $G=T^l$ with $l\geqslant2$, and let 
$$R=\{(t^x,1,\dots,1),(1,t^x,\dots,1),\dots,(1,1,\dots,t^x),
(t^x,t^x,\dots,t^x)\mid x\in T\}.$$
Then $\Cay(G,R\cup R^{-1})$ is subnormal and half-symmetric.
\end{theorem}

Most finite simple groups $T$ contain elements which are not conjugate
in $\Aut(T)$ to their inverses.
Here is an example.
Let $T=\Sz(q)$ with $q=2^{2e+1}\geqslant8$, and let $t$ be an element 
of $T$ of order 4.
Then $t$ is not conjugate in $\Aut(T)$ to $t^{-1}$.

\section{Proof of Theorem~\ref{subnormal-2-trans}}

Let $\Ga=(V,E)$ be a digraph.
For $v\in V$, let 
$\Ga(v)=\{w\in V\mid (v,w)\ \mbox{is an arc of $\Ga$}\}$.
Let $G_v^{[1]}$ be the kernel of $G_v$ acting on $\Ga(v)$.
Then $G_v^{[1]}$ is normal in $G_v$. Let $\Ga_{0,1,\cdots,i}(v)=\{w\mid \mbox{the distance between }v \mbox{ and } d\mbox{ are not larger than }i\mbox{ in }\Ga\}$. 
We first prove a simple lemma about vertex stabilisers.

\begin{lemma}\label{vt-stabilizer}
Let $\Ga$ be a connected $G$-vertex transitive digraph.
Then for a vertex $v$ and a normal subgroup $N\lhd G$, if $N_v^{\Ga(v)}$ is semiregular, then 
$N_v\cong N_v^{\Ga(v)}$ is faithful.
\end{lemma}
\proof 
Suppose that $N_v^{\Ga(v)}$ is semiregular for a $v$. Since $N\lhd G$, for any $w\in \Ga$, there is an element $g\in G$ such that $w=v^g$. Thus  $N_{v^g}=N_v^g$ and $N_w^{\Ga(w)}$ is semiregular. For the contrary, suppose that there exists an  $x\in N$ such that $x$ fixes pointwesely $\Ga_{0,1}(v)$. Let $i\ge 1$ be the maximal integer such that $x$ fixes  $\G_{0,1, \cdots , i}(v)$  but moves a vertex $w'\in \Ga_{i+1}(v)$. Let $v'\in \Ga_{i-1}(v)$, $w\in \Ga_i(v)$ and $w'\in\Ga_{i+1}(v)$ such that $(v',w,w')$ is a 2-arc.  Then $x$ fixes $v'$, $w$, and  moves $w'$. Thus $x\in G_{wv'}$ and acts non-trivially on $\Ga_1(w)$. So  $N_w^{\Ga(w)}$ is not semiregular, a contradiction.\qed

For a group $X$ and a core free subgroup $H\leqslant X$, denote by $[X:H]$
the set of right cosets of $H$ in $X$, that is
$$[X:H]=\{Hx\mid x\in X\}.$$
For any subset $S\subset X$, define the {\it coset graph} of $X$ with 
respect to $H$ and $S$ to be the digraph $\Ga$ with vertex set $[X:H]$
and such that two vertices $Hx,Hy\in V$ are adjacent, written as
$Hx\sim Hy$, if and only if $yx^{-1}\in HSH$;
denoted by $\Ga=\Cos(X,H,HSH)$.
Then $X\leqslant\Aut\Ga$, and $\Ga$ is $X$-vertex transitive.
For convenience, write $H\{g\}H=HgH$, where $g\in X$.
The following properties are known and easy to prove.

\begin{lemma}\label{coset-graphs}
Let $X$ be a group, $H$ a core free subgroup, and $g\in X$.
\begin{itemize}

\item[(i)] $\Cos(X,H,HgH)$ is connected if and only if $\l H,g\r=X;$
\item[(ii)] $\Cos(X,H,HgH)$ is $X$-edge transitive;

\item[(iii)] $\Ga=\Cos(X,H,H\{g,g^{-1}\}H)$ is undirected and $X$-edge
 transitive;
further, $\Ga$ is $X$-arc transitive if and only if $HgH=Hg^{-1}H$.
\end{itemize}
\end{lemma}

Let $\Aut(X,H)=\l\s\in\Aut(X)\mid H^\s=H\r$.
Then an element $\s\in\Aut(X,H)$ acts on $[X:H]$ by $(Hx)^\s=Hx^\s$.
Let $\s\in\Aut(X,H)$ be such that $(HgH)^\s=HgH$.
Then for any two vertices $Hx,Hy$, we have
$$\begin{array}{rcl}
Hx\sim Hy & \Leftrightarrow & yx^{-1}\in HgH\\
          & \Leftrightarrow 
   & y^\s (x^\s)^{-1}=(yx^{-1})^\s\in (HgH)^\s=HgH\\
          & \Leftrightarrow & Hx^\s\sim Hy^\s\\
\end{array}$$
Thus $\s$ maps all edges to edges, and so $\s$ induces an
automorphism of $\Ga$.

\begin{lemma}\label{Aut(X,H,HgH)}
Let $\Ga=\Cos(X,H,HgH)$, and $\s\in\Aut(X,H)$.
If $(HgH)^\s=HgH$, then $\s$ induces an automorphism of $\Ga$.
\end{lemma}

For a group $G$, the symmetric group $\Sym(G)$ contains two regular subgroups $\widehat G$ and $\check G$, where
$$\check G=\{\check g:\ x\mapsto g^{-1}x\ \mbox{for all $x\in G$}
             \mid g\in G\},$$
consisting of left multiplications of elements $g\in G$ and $\hat G$ with 
$$\hat G=\{\hat g:\ x\mapsto x g\ \mbox{for all $x\in G$}
             \mid g\in G\},$$ consisting of left multiplications of elements $g\in G$. 
Then $\N_{\Sym(G)}(\widehat G)=\widehat G\rtimes\Aut(G)$, the {\it holomorph} of $G$, and 
$\widehat G\C_{\Sym(G)}(\widehat G)=\widehat G\circ\check G=\widehat G\rtimes\Inn(G)$.

For a subset $S\subset G$, let
$$\Aut(G,S)=\{\s\in\Aut(G)\mid S^\s=S\}.$$
Then $\Aut(G,S)\leqslant\Aut(G)\leqslant\Sym(G)$, and as subgroups of $\Sym(G)$, 
it is easily shown that $\Aut(G,S)$ normalizes $\widehat G$.
Moreover, for the Cayley graph $\Ga=\Cay(G,S)$, 
by \cite[Lemma~2.1]{Godsil}, we have
$$\N_{\Aut\Ga}(\widehat G)=\widehat G\rtimes \Aut(G,S).$$

The subgroup $\Aut(G,S)$ plays an important role in the study of Cayley graphs.  Assume that If $\widehat G\lhd X\leqslant Aut\Ga$. Then $X_\a\leqslant \Aut(G,S)$ where $\a$ is a vertex of $\Ga$. A special type of normal Cayley graph satisfies $X_\a\geqslant \Inn(G,S)$, in this case, we call $\Ga$ a {\it holomorph Cayley graph.}

Suppose $\Ga=\Cay(G,S)$ is a holomorph graph with $H=\widehat G\circ\check G=\widehat G\rtimes\Inn(G)$.
Let $\b\in \Ga(\a)=S$, let $g\in H_{\a\b}$, then $\b^h=\b$, that is $h\in \C_{G}(\b)$.  On the contrary, if $h\in \C_{ G}(\b)$, then $\b^h=\b$, so $h\in  H_{\a\b}$. Thus $H_{\a\b}=\C_{ G}(\b).$ Thus the following lemma holds.

\begin{lemma}\label{homo}
Suppose $\Ga=\Cay(G,S)$ is a holomorph with $H=\widehat G\circ\check G.$ Then $H_{\a\b}=\C_{ G}(\b).$
\end{lemma}

The next lemma shows that, for a prime $p$ and an integer $d$, 
a complete bipartite graph $\bfK_{p^d,p^d}$ is a $2$-arc transitive subnormal Cayley graph.

\begin{lemma}\label{example}
Let $\Ga=\bfK_{p^d,p^d}$, where $p$ is an odd prime and $d\geqslant1$.
Then $\Ga\cong\Cay(G,S)$, where $G\cong\ZZ_p^d\rtimes\ZZ_2$ and $S$
consists of all involutions of $G$, and there exist subgroups $X,Y<\Aut\Ga$ such that 
$\widehat G\lhd X\lhd Y<\Aut\Ga$, $X=\widehat G\rtimes\Aut(G)$, and $Y/X\cong\ZZ_2$.
\end{lemma}
\proof
Let $G=N\rtimes\l z\r\cong\ZZ_p^d\rtimes\ZZ_2$, where $p$ is an odd
prime and $z$ reverses every element of $N$, that is, for each element
$x\in N$, $x^z=x^{-1}$.
Let $S=G\setminus N$, and let $\Ga=\Cay(G,S)$.
Then $S$ consists of all involutions of $G$. Let $V_1$ be the vertex set corresponding to the elements in $N$, $V_2$ be the vertex set corresponding to the elements in $G\setminus N$. Then each vertex in $V_1$ is adjcent to all vertices in $V_2$ and each vertex in $V_2$ is adjcent to all vertices in $V_1$ as well. So 
$\Ga\cong\bfK_{p^d,p^d}$.
Thus, $\Aut\Ga\cong\S_{p^d}\wr\S_2$.
Further, $\Aut(G,S)=\Aut(G)\cong\AGL(d,p)=\ZZ_p^d\rtimes\GL(d,p)$, and
$\Aut(G,S)$ acts 2-transitively on $S$. 
 
Let $X=\N_{\Aut\Ga}(\widehat G)$, and 
let $C=\widehat G\C_{\Aut\Ga}(\widehat G)$.
Then $X=\widehat G\rtimes\Aut(G,S)=\widehat G\rtimes\Aut(G)$, and 
$C=\widehat G\times\check G$.
Thus $\Ga$ is $(X,2)$-arc transitive and $C$-arc transitive.
Let $v$ be the vertex of $\Ga$ corresponding to the identity of $G$.
Then $C_v=\{(\widehat g,\check g)\mid g\in G\}\cong G$. Let 
$\Ga'=\Cos(C,C_v,C_v(\widehat z,1)C_v)$ and $\phi$ a map from vetices of $\Ga$ to vertices of $\Ga'$ such that for any vertex $C_vx\in V\Ga'$ and $x\in V\Ga$, $\phi: C_vx\mapsto x$. Then $\phi$ is an isomorphism of $\Ga$ to $\Ga'$. Thus $\Ga\cong\Cos(C,C_v,C_v(\widehat z,1)C_v).$

We label $\Aut(\widehat G)=\{\widehat x\mid x\in\Aut(G)\}$, and
$\Aut(\check G)=\{\check x\mid x\in\Aut(G)\}$.
Then $\Aut(C)=\Aut(\widehat G\times\check G)
=(\Aut(\widehat G)\times\Aut(\check G)).\l\t\r$, where 
$\t:\ (\widehat x,\check y)\mapsto(\widehat y,\check x)$ for all
$(\widehat x,\check y)\in\Aut(\widehat G)\times\Aut(\check G)$.
Let  $(\widehat x,\check y)\in\Aut(C)$ normalize
$C_v=\{(\widehat g,\check g)\mid g\in G\}$.
Then $(\widehat g^{\widehat x},\check g^{\check y})\in C_v$ for any $g\in G$.
 Thus $g^{yx^{-1}}=g$ for any $g\in G$, that is $yx^{-1}\in \mathbb Z(G)=1 .$
 Hence $x=y$ and
$\Aut(C,C_v)=\l(\widehat x,\check x)\mid x\in\Aut(G)\r\times\l\t\r$.
Since $C\lhd X$ and $\C_X(C)=1$, it follows that $X\leqslant\Aut(C)$.
Further, $C_v\lhd X_v\leqslant\Aut(C,C_v)$, and it follows that $\Aut(C,C_v)=X_v\times\l\t\r$.
Noticing that $(\widehat z,\check z)\in C_v$ and $\check z$ is an involution, we have 
$$(C_v(\widehat z,1)C_v)^\t=C_v(\widehat z,1)^\t C_v=C_v(1,\check z)C_v
   =C_v(1,\check z)(\widehat z,\check z)C_v=C_v(\widehat z,1)C_v.$$
By Lemma~\ref{Aut(X,H,HgH)}, $\t\in\Aut\Ga$ and $\Aut(C,C_v)<\Aut\Ga$.
Now $Y:=C\Aut(C,C_v)$ is such that $|Y:X|=2$.
We obtain that $\widehat G\lhd X\lhd Y<\Aut\Ga$.
Since $\t\in Y$ does not normalizes $\widehat G$, $\widehat G$ is not normal in $Y$.

Therefore, as $\Ga$ is $(Y,2)$-arc transitive, 
$\bfK_{p^d,p^d}$ is a $2$-arc transitive subnormal Cayley graph.
\qed

The following is a property regarding $2$-transitive permutation groups,
which is obtained by inspection of the classification of 2-transitive
permutation groups, refer to \cite{DM-book}. 

\begin{lemma}\label{HA-stab}
Let $X$ be a $2$-transitive permutation group on $\Ome$. 
Then the socle of X is either a regular elementary abelian $p$-group, or a nonregualr nonabelian simple group.

Furter, assume that $N\lhd\lhd X$ is imprimitive on $\Ome$.
Then $X$ is affine with $\soc(X)=\ZZ_p^e$, where $p$ is a prime and
$e\geqslant1$, and further, the following hold:
\begin{enumerate}
\item[(i)] Either $N\leqslant\soc(X)$, or $\ZZ_p^e.\ZZ_b\cong N\lhd X$ and 
$N$ is a Frobenius group, where $b$ divides $p^{e'}-1$ and $e'$ is a
  proper divisor of $e$. 

\item[(ii)] $X_\o$ has no non-trivial normal subgroup of $p$-power
  order for $\o\in\Ome$.
\end{enumerate}
\end{lemma}

This has an application to 2-arc transitive graphs.

\begin{lemma}\label{normal-subgp-2-trans}
Let $\Ga$ be a $(Y,2)$-arc transitive graph, and let $H$ be a
subnormal subgroup of $Y$ which is vertex transitive on $\Ga$.
Then either $H_v^{\Ga(v)}$ is center free and $\Ga$ is $H$-arc transitive, or 
$H_v$ is abelian and acts faithfully and semiregularly on $\Ga(v)$.
\end{lemma}
\proof
Since $H\lhd\lhd Y$, we have that $H_v\lhd\lhd Y_v$, and 
$H_v^{\Ga(v)}\lhd\lhd Y_v^{\Ga(v)}$ and $Y_v^{\Ga(v)}$ is a 2-transitive permutation group.
If $H_v^{\Ga(v)}$ is primitive, then $\Ga$ is $H$-arc transitive and $H_v^{\Ga(v)}\geqslant \soc(Y_v^{\Ga(v)})$ is center free by Lemma~\ref{HA-stab}.

Now suppose $H_v^{\Ga(v)}$ is imprimitive. Since $Y_v^{\Ga(v)}$ is a 2-transitive permutation group, it follows
from Lemma~\ref{HA-stab} that either $H_v^{\Ga(v)}\leqslant\soc(Y_v^{\Ga(v)})\cong\ZZ_p^e$, where $p$ is a prime
and $e\geqslant1$ or
$\soc(Y_v^{\Ga(v)})=\ZZ_p^e\leqslant H_v^{\Ga(v)}=\ZZ_p^e.\ZZ_b$ and $H_v^{\Ga(v)}$ is center free.

For the former, since $\soc(Y_v^{\Ga(v)})$ is regular, $H_v^{\Ga(v)}$
is semiregular.   
By Theorem~\ref{vt-stabilizer}, $H_v\cong H_v^{\Ga(v)}$ is faithful and
abelian.

For the latter, since $Y_v^{\Ga(v)}$ is 2-transitive, we have that
$H_v^{\Ga(v)}\geqslant\soc(Y_v^{\Ga(v)})$ is transitive, and hence $\Ga$ is
$H$-arc transitive.
\qed

To prove Theorem~\ref{subnormal-2-trans}, we need the next property on permutation groups.
 
\begin{lemma}\label{G1=G2}
Let $G_1,G_2<\Sym(\Ome)$ be regular which normalizes each other.
If $G_1/(G_1\cap G_2)$ is abelian, then $G_1=G_2$.
\end{lemma}
\proof
Let $X=G_1G_2$, and $C=G_1\cap G_2$.
Then $C$ is semiregular on $\Ome$, and $C\lhd X$.
Let $\ov G_1=G_1/C$, $\ov G_2=G_2/C$, and $\ov X=X/C$.
Let $\Ome_C$ be the set of $C$-orbits on $\Ome$.
Then both $\ov G_1$ and $\ov G_2$ are regular on $\Ome_C$
as $G_1,G_2$ are both regular on $\Ome$.

Suppose that $G_1\not= G_2$.
Then $\ov G_i\not=1$, and $\ov X=\ov G_1\times\ov G_2$.
In particular, $\ov G_2\leqslant\C_{\Sym(\Ome_C)}(\ov G_1)$.
If $\ov G_1$ is abelian, then 
$\ov G_2\leqslant\C_{\Sym(\Ome_C)}(\ov G_1)=\ov G_1$.
Thus $\ov G_2=\ov G_1$, and so $G_1=G_2$, which is a contradiction.
\qed

Now we are ready to prove Theorem~\ref{subnormal-2-trans}.

\vskip0.1in
\noindent{\bf Proof of Theorem~\ref{subnormal-2-trans}:}
Let $\Ga=\Cay(G,S)$ be a $(Y,2)$-arc transitive graph with vertex set $V$. Then $G$ is regular on $V$.

If $G\lhd Y$ then the theorem holds.
Now we supppose $G\ntriangleleft Y$. Then $\N_{\Aut\Ga}(G)<Y.$ Let $X$ be the maximal subnormal subgroup of  $Y$ contained in $\N_{\Aut\Ga}(G)$, we have $G\lhd X\lhd\lhd Y$.  If $X\lhd Y$ then $\N_Y(X)=Y>\N_{\Aut\Ga}(G)$, otherwise there is a group $K>X$ such that  $X\lhd K\lhd\lhd Y$, so  $\N_Y(X)\geq K$ with $K\cap \N_{\Aut\Ga}(G)=X$ as $X$ is maximal. Thus $\N_Y(X)\neq \N_{\Aut\Ga}(G).$
Since $G<X<\N_Y(X)\leqslant Y$, any element $y\in\N_Y(X)\setminus\N_Y(G)$ is
such that $G^y\not=G$ and $X^y=X$.
 
Let $C=G\cap G^y$.
Then for any $x\in X$, we have $G^x=G, (G^y)^x=G^{yx}=G^{x'y}=G^y$ for some $x'\in X$. Thus  $C^x=C$ and $C,G, G^y\lhd X$; in particular, $G$ and $G^y$ normalizes each other.
Let $\ov G=G/C$ and $\ov G^y=G^y/C$, let $V_C$ be the set of $C$-orbits on $V$.
By Lemma~\ref{G1=G2}, $\ov G$ is not abelian as $G\not= G^y.$
Since $G,G^y$ are both regular on $V$, the subgroup $C$ is semiregular
on $V$, and $\ov G$, $\ov G^y$ are both regular on $V_C$.
Further, $\ov G^y\leqslant\C_{\Sym(V_C)}(\ov G)$.

Let $H=GG^y$.
Then $H\lhd X$.
Let $\ov H=H/C$, and $\ov X=X/C$.
Then $\ov G\times \ov G^y=\ov H\lhd \ov X$.
Let $v$ be a vertex of $\Ga$.
Then $H=G{:}H_v=G^y{:}H_v$, and $H_v\cong H/G^y\cong G/C=\ov G$.
Further, since $G<H\lhd X\lhd\lhd Y$, we have 
$1\not=H_v\lhd X_v\lhd\lhd Y_v$, and 
$1\not=H_v^{\Ga(v)}\lhd X_v^{\Ga(v)}\lhd\lhd Y_v^{\Ga(v)}$.
By Lemma~\ref{HA-stab}, we conclude that either
$\soc(Y_v^{\Ga(v)})\leqslant H_v^{\Ga(v)}$, or $Y_v^{\Ga(v)}$ is 
affine with socle isomorphic to $\ZZ_p^d$,
$H_v^{\Ga(v)}<\soc(Y_v^{\Ga(v)})\cong\ZZ_p^d$, and
$H_v^{\Ga(v)}$ is semiregular.

Let $\a$ be the vertex of $\Ga_C$ containing $v$, that is, $\a=v^C$.
Then the stabilizer $\ov H_{\a}$ is isomorphic to $\ov G$ 
as $\ov G\times\ov G^y=\ov H=\ov G{:}\ov H_{\a}=\ov G^y{:}\ov H_{\a}$.
On the other hand, $\ov H_{\a}$ is isomorphic to a factor group of 
$H_v$, that is, $\ov H_{\a}\cong H_vC/C\cong H_v/(H_v\cap C)$.

Suppose that $H_v$ is abelian.
Then the factor group $\ov H_{\a}\cong H_v/(H_v\cap C)$ is abelian.
Since $\ov G\cong\ov G^y\cong\ov H_{\a}$, we conclude that $G$ is
abelian by Lemma~\ref{G1=G2}, which is a contradiction.
Thus, $H_v$ is not abelian.
By Lemmas~\ref{HA-stab} and ~\ref{normal-subgp-2-trans}, either $Y_v^{\Ga(v)}$ is almost
simple, or $H_v^{\Ga(v)}=\ZZ_p^d{:}H_o$ is a Frobenius group.
In particular, $H_v$ is transitive on $\Ga(v)$, and $\Ga$ is $H$-arc
transitive.
 
 Since $\ov G\times \ov G^y=\ov H\leqslant \Aut(\Ga_C)$, and $\ov G$ is regular on $\Ga_C$, we have $\Ga_C$ is a holomorph graph $\Cay(\ov G, S)$. So $ \ov H_\a^{\Ga_C(\a)}\cong\ov H_\a=\Inn(\ov G,S).$
 Suppose that $Y_v^{\Ga(v)}$ is almost simple.
Then $\soc(Y_v^{\Ga(v)})\leqslant H_v^{\Ga(v)}\leqslant Y_v^{\Ga(v)}$.
Since $Y_v^{\Ga(v)}$ is 2-transitive, 
either $H_v^{\Ga(v)}$ is 2-transitive, or $|\Ga(v)|=28$,
$H_v^{\Ga(v)}\cong\PSL(2,8)$ and $Y_v^{\Ga(v)}\cong\PGammaL(2,8)$.
For the former, the graph $\Ga_C$ is a holomorph 2-arc transitive
graph, which is not possible, see \cite[Theorem\,1.3]{ET-Cay}.
For the latter, since $|\Ga_C(\a)|=|\Ga(v)|=28$, we have $\ov H_\a=\D_{18}$ which have index 28 in $\ov H_\a=\PSL(2,8)$. However $\ov H_\a=\D_{18}$ is not the centraliser of any element in $\ov H_\a=\PSL(2,8)$,  which is not possible.
 
Thus, $\Ga$ and $\Ga_C$ are of valency $p^d$, and 
$\ov G\cong H_v^{\Ga(v)}=\ZZ_p^d{:}H_o\cong\ZZ_p^d{:}\ZZ_b$ is a Frobenius group;
in particular, $\ov G$ is center free.
Hence $\ov H\cong(\ZZ_p^d{:}H_o)\times(\ZZ_p^d{:}H_o)$.
Now $\Ga_C$ is a holomorph Cayley graph of $\ov G=\ZZ_p^d{:}H_o$, 
that is, $\Ga_C=\Cay(\ov G,S)$ such that $S$ is a full conjugacy class of elements of $\ov G$, and $|S|=p^d$.
Let $\a$ to be the vertex of $\Ga_C$ corresponding to the identity of 
$\ov G$, and let $\b\in\Ga_C(\a)=S$. 
Then $\ov H_\a\cong\ZZ_p^d{:}H_o$, and $\ov H_{\a\b}\cong H_o$.
Since $\Ga_C$ is undirected, we have $S=S^{-1}$ and so $\b^{-1}\in S$. 
Now $\C_{\ov G}(\b)\cong\ov H_{\a\b}$ so $\b$ is not order $p$. Further as $\ov G=\ZZ_p^d{:}\ZZ_b$ is a Frobenius group, $\b$ is not conjugate to $\b^{-1}$ if $o(\b)>2$. Hence $\b$ is an involution.
It follows that $p$ is odd and $H_o=\l\b\r\cong\ZZ_2$.
So $\Ga_C\cong\bfK_{p^d,p^d}$, and $\Ga$ is a normal cover of $\Ga_C$. By Lemma~\ref{example}, the theorem holds.
\qed

\section{Subnormal transitive subgroups}

Let $G\lhd\lhd X\leqslant\Sym(\Ome)$ be such that $G$ is transitive on $\Ome$.
Assume that $G=N_0\lhd N_1\lhd\dots\lhd N_r=X$,
where $N_{i+1}=\N_X(N_i)>N_i$.
A natural question is whether $r$ has an upper-bound.
For characteristic simple groups, we have a positive answer.

\begin{lemma}\label{sub-normality}
Let $G\leqslant\Sym(\Ome)$ be a finite characteristic simple group.
If $G\lhd\lhd X\leqslant\Sym(G)$ and $G$ is transitive, then either $G\lhd X$, or
there exists a group $N$ such that $G\lhd N\lhd X$.
\end{lemma}
\proof
Write $G=T^k$, where $T$ is a simple group and $k\geqslant1$.
Suppose that $G\lhd\lhd X\leqslant\Sym(G)$ and $G$ is not normal in $X$.
Let $N=\N_X(G)$.
Then $N<X$, and there exists $x\in X\setminus N$ such that 
$N^x=N$ and $G^x\not=G$.
Let $C=G\cap G^x$ and $H=GG^x$. Then $C,G,$ and $G^x$ are normal in $N$, in particular, $G$ and $G^x$ normalizes each other.
If $G$ is abelian, then $G$ is regular and $G/C$ is abelian,
which is a contradiction to Lemma~\ref{G1=G2} since now $G^x\not=G$.
Thus $G$ and so $T$ is nonabelian.

Let $G=N_0\lhd N_1\lhd N_2\lhd\dots\lhd N_r=X$.
Let $M_i=\l G^x\mid x\in N_i\r$, where $2\leqslant i\leqslant r$.
We claim that $M_i=G\times T^{m_i}$ for some positive integer $m_i$.
First, $M_2=\l G^x\mid x\in N_2\r\lhd N_2$.
Since $G^x\lhd N_1$ for $x\in N_2$, we conclude that $GG^x=T^n$ 
for some $n>k$, and as $G\lhd GG^x$, we have $GG^x=G\times T^l$.
It follows that $M_2=G\times T^{m_2}$ for some positive integer $m_2$.
Assume inductively that $M_i=G\times T^{m_i}$ for some positive
integer $m_i$. 
Then $M_i=T^{k+m_i}$ is a characteristically simple group.
Arguing as for $M_2$, with $M_i$ in the position of $G$, we obtain 
$M_{i+1}=\l M_i^x\mid x\in N_{i+1}\r=M_i\times T^n=G\times T^{m_{i+1}}$, 
where $m_{i+1}$ is a positive integer.
By induction, $M_r=G\times T^{m_r}$, and hence $G\lhd M_r\lhd N_r=X$.
\qed

However, we have been unable to extend this lemma for general groups.

\begin{question}\label{index-subnormal}
{\rm
Let $G=N_0\lhd N_1\lhd\dots\lhd N_r=X\leqslant\Sym(\Ome)$, 
where $N_{i+1}=\N_X(N_i)>N_i$.
Assume that $G$ is transitive.
Is it true that $r\leqslant2$?
}
\end{question}

In the rest of this section, we construct a family of half-symmetric
graphs which are subnormal Cayley graphs, and prove
Theorem~\ref{half-symmetric}. 

Let $T$ be a nonabelian simple group, and let $k\geqslant2$.
Let 
$$X=T^k.(\Out(T)\times\S_k)$$ 
be a primitive permutation group on
$\Ome\equiv T^{k-1}$ of simple diagonal type, see \cite{LPS-alt}.
Then the stabilizer
$$X_\o=D.(\Out(T)\times\S_k)=D.\Out(T)\times\S_k,$$ 
where $D.\Out(T)=\{(t,t,\dots,t)\mid t\in \Aut(T)\}$, 
and the socle $M:=\soc(X)=T^k=T_1\times T_2\times\dots\times T_k$.
Let $G\times\{1\}=T_1\times\dots\times T_{k-1}\times\{1\}\lhd\soc(X)$, and $N=\N_X(G\times\{1\})$.
Then $G\times\{1\}$ is regular on $\Ome$, and $N=T^k.(\Out(T)\times\S_{k-1})$.

\vskip0.1in
\noindent{\bf Proof of Theorem~\ref{half-symmetric}:}
Using the notation defined above, assume further that $k\geqslant3$, and 
there exists an element $t\in T$ such that $t$ is not
conjugate in $\Aut(T)$ to $t^{-1}$.
Let $g=(t,1,\dots,1,1)\in G\times\{1\}$ where $t\in T$, and let
$$\Ga=\Cos(X,X_\o,X_\o\{g,g^{-1}\}X_\o).$$ 
Then $\Ga$ is $X$-edge transitive, and $X\leqslant\Aut\Ga\leqslant\Sym(\Ome)$.
Further, $\Ga$ is not a complete graph, and so $\Aut\Ga\not=\Sym(\Ome)$.
By \cite{LPS-alt}, $\Aut\Ga=X$.

Suppose that $\Ga$ is arc-transitive.
Then by Lemma~\ref{coset-graphs}, $X_\o gX_\o=X_\o g^{-1}X_\o$, and 
so $g=xg^{-1}y$, for some elements $x,y\in X_\o$.
Since $X_\o=D.\Out(T)\times\S_k$, the elements 
$x=(t_1,t_1,\dots,t_1)\pi_1$, and $y=(t_2,t_2,\dots,t_2)\pi_2$, 
where $t_i\in\Aut(T)$, and $\pi_i\in\S_k$.
Thus 
$$(t,1,\dots,1)=g=xg^{-1}y =
(t_1,t_1,\dots,t_1).\pi_1(t^{-1}t_2,t_2,\dots,t_2)\pi_1^{-1}.\pi_1\pi_2.$$
It follows that $\pi_1\pi_2=1$, and the element on the right hand side has
exactly one entry equal to $t_1t^{-1}t_2$ and the other entries equal
to $t_1t_2$.
Since $k\geqslant3$, we conclude that $t_1t_2=1$ and $t=t_1t^{-1}t_2$.
Thus $t=t_2^{-1}t^{-1}t_2$ and so $t$ is conjugate to $t^{-1}$ which is a contradiction.
Hence $\Ga$ is half-symmetric.

Finally, since $G\times\{1\}$ is regular on $\Ome$, 
$\Ga$ is a Cayley graph of $G\times\{1\}$,
that is, $\Ga=\Cay(G\times\{1\},S\times\{1\})$ for some subset 
$S\times\{1\}\subset G\times\{1\}$.
Let $\o$ be the vertex corresponding to $X_\o$, let $\b=X_\o g$. Then the stabilizer of $\o$ in $X$ is $X_\o$, and the stabilizer of $\b=X_\o^g.$ So $X_{\o\b}=X_\o\cap X_\o^g=\C_{\Aut(T)}(g)\times\S_{k-1}$.
Since $g$ is not conjugate to $g^{-1}$, we have $\Ga$ is not $X$-arc transitive. By Lemma 2.1 in \cite{Li-Song}, $X_\o$ have two orbits of the same size on $\Ga(\o)$, and each have size
$\val(\Ga)=|X_\o:X_{\o\b}|=|\Aut(T):\C_{\Aut(T)}(g)|.k.$ Thus $$\val(\Ga)=2|X_\o:X_{\o\b}|=2|\Aut(T):\C_{\Aut(T)}(g)|.k.$$
Let $\pi=(12\dots k)$, and let $g_i=g^{\pi^i}$.
Then the $i$-th entry of $g_i$ is $t$ and the others equal 1, and 
$$X_\o\{g,g^{-1}\}X_\o=\{X_\o g_i^x, X_\o(g_i^{-1})^x\mid x\in X_\o\}.$$
Note that $g_i^x=(1,\dots,t^{x_i},\dots,1)$, where $x_i\in\Aut(T_i)$.
For $i\leqslant l=k-1$, let $\ov g_i^x$ be the projection of $g_i^x$ in 
$G=T_1\times\dots\times T_{k-1}$.
For $i=k$, $g_k^x$ has the following property
$$g_k^x=(1,\dots,1,t^{x_k})\equiv
   ((t^{x_k})^{-1},\dots,(t^{x_k})^{-1},1)\ (\mod X_\o).$$
Let $\ov g_k^x=((t^{x_k})^{-1},\dots,(t^{x_k})^{-1})$ be the projection
of $g_k^x$ in $G$.
Then all $\ov g_i^x$ for $1\leqslant i\leqslant k$ lie in $S$.
Similarly, we have the projections $(\ov g_i^{-1})^x$ of
$(g_i^{-1})^x$.
Since $|S|=|S\times\{1\}|=\val(\Ga)=2|\Aut(T):\C_{\Aut(T)}(g)|.k$, it follows that 
$$S=\{\ov g_i^x, (\ov g_i^{-1})^x \mid 1\leqslant i\leqslant k,\ x\in X_\o\}.$$
Thus $\Ga$ can be represented as a Cayley graph of $G$, that is,
$\Ga\cong\Cay(G,S)\cong\Cay(G\times\{1\},S\times\{1\})$.
As $G\cong G\times\{1\}\lhd\soc(X)\lhd X$, the Cayley graph $\Ga$ is subnormal and has the form stated in Theorem~\ref{half-symmetric}.
\qed

\end{document}